\documentclass[10pt,a4paper]{amsart}
\DeclareSymbolFont{AMSb}{U}{msb}{m}{n}
\DeclareSymbolFontAlphabet{\Bbb}{AMSb}
\usepackage{latexsym}
\usepackage[dvips]{graphicx}
\usepackage{amstext,amsfonts,amsmath,amsthm,graphicx,amssymb,amscd,epsfig}
\usepackage{psfrag} % coloca texto TeX dentro de las figuras necesita los
\input epsf
\newcommand{\C}{\mathbb{C}}
\newcommand{\R}{\mathbb{R}}

\newcommand{\newc}{\newcommand}

\def\be{\begin{equation}}
\def\ee{\end{equation}}
\def\bq{\begin{eqnarray}}
\def\eq{\end{eqnarray}}
\def\beq{\begin{eqnarray}}
\def\eeq{\end{eqnarray}}
\def\ba{\begin{array}}
\def\ea{\end{array}}
\def\bi{\begin{itemize}}
\def\ei{\end{itemize}}

\newc{\e}{{\mbox{e}}}
\newc{\deter}{\operatorname{{det}}}  \newc{\const}{\operatorname{const}}
\newc{\diam}{\operatorname{{diam}}}  \newc{\supp}{\operatorname{{supp}}}
\newc{\grad}{\operatorname{{grad}}}  \newc{\dist}{\operatorname{{dist}}}
\newc{\graph}{\operatorname{{graph}}}\newc{\length}{\operatorname{{length}}}
\newc{\id}{\operatorname{{id}}}      \newc{\card}{\operatorname{{card}}}

\newcommand{\spec}{\operatorname{{Spec}}}

\begin{document}
\title[A solution to the Markus-Yamabe Conjecture]
{The Markus-Yamabe Conjecture for differentiable  vector fields of
$\R^2$. }

\author{Carlos Gutierrez and  Roland Rabanal}

%\date{\today}

\address{Carlos Gutierrez and Roland Rabanal \newline
Instituto de Ci\^encias Matem\'aticas e de Computa\c{c}\~ao,
Universidade de S\~ao Paulo, caixa postal 668, 13560-970, S\~ao
Carlos SP, Brasil} \email{gutp@icmc.usp.br, roland@icmc.usp.br}

\dedicatory{Dedicated to C\'esar Camacho on his 60th Birthday}

%\keywords{Planar Vector Fields, Global Injectivity, Asymptotic Stability}%

\begin{abstract}

(a) Let $X\colon \R^2 \rightarrow \R^2$ be a differentiable map
(not necessarily $C^1$) and let $\spec(X)$ be the set of (complex)
eigenvalues of the derivative $DX_p$ when $p$ varies in $\R^2$.
If, for some $\epsilon>0$, $\spec(X)\cap [0,\epsilon)=\emptyset$
then  $X$ is injective.

(b) Let $X\colon \R^2 \rightarrow \R^2$ be a differentiable vector
field such that $X(0)=0$ and $\spec(X)\subset \{z \in \C : \Re(z)
< 0\}.$ Then, for all $p\in \R^2,$ there is a unique positive
trajectory starting at $p;$ moreover the $\omega-$limit set of $p$
is equal to $\{0\}.$
\end{abstract}

\maketitle

%%%%%%%%%%%%%%%%%%%%%%%%%%%%%%%%%%%%%%%%%%%%%%%%%%%%%%%%%%%%%%%%

%%%%%%%%%%%%%%%%%%%%%%%%%%%%%%%%%%%%%%%%%%%%%%%%%%%%%%%%%%%%%%%%

\end{document}